\documentclass[11pt]{article}
\usepackage{latexsym,amssymb}
\def\cE{{\cal E}}
\def\cF{{\cal F}}
\def\cG{{\cal G}}
\def\cR{{\cal R}}
\def\cS{{\cal S}}
\def\cL{{\cal L}}
\def\cW{{\cal W}}
\def\cO{{\cal O}}
\def\CC{\mathbb C}
\def\LG{\mathbb {LG}}
\def\LF{{\mathbb {LF}}}
\def\ZZ{\mathbb Z}
\def\GG{\mathbb G}
\def\PP{\mathbb P}
\def\QQ{\mathbb Q}
\def\FF{\mathbb F}
\def\l{\lambda}
\def\a{\alpha}
\def\m{\mu}
\def\e{\varepsilon}
\def\ot{{\mathord{\otimes }\,}}
\def\op{{\mathord{\oplus }\,}}
\def\otc{{\mathord{\otimes\cdots\otimes}\;}}
\def\pc{{\mathord{+\cdots +}}}
\def\lra{{\mathord{\;\longrightarrow\;}}}
\def\ra{{\mathord{\;\rightarrow\;}}}
\def\we{{\mathord{\wedge}}}

\newcommand\qed{{\hspace*{\fill} $\Box$}} 
\newtheorem{theo}{Theorem}[section]
\newtheorem{coro}[theo]{Corollary}
\newtheorem{lemm}[theo]{Lemma}
\newtheorem{prop}[theo]{Proposition}

\begin{document}

\title{On branched coverings of some  homogeneous spaces}
\author{Meeyoung Kim, Laurent Manivel}
\date{}
\maketitle

\centerline{\it Dedicated to the memory of Prof. Michael Schneider}

\begin{quote}{\small {\sc Abstract}. {\em
We study nonsingular branched coverings of a homogeneous space $X$. 
There is a vector bundle associated with
such a covering which was conjectured by O. Debarre to be ample when the
Picard number of $X$ is one. We prove this conjecture, which implies 
Barth-Lefschetz type theorems, for lagrangian grassmannians, and for 
quadrics up to dimension six. We propose a conjectural extension to
homogeneous spaces of Picard number larger than one and prove a
weaker version.}}\end{quote}

\section{Introduction and main results}

Let $f: Y\ra\PP^n$ be a branched covering of degree $d$ of a complex
projective space where $Y$ is a nonsingular connected complex
projective variety. A
celebrated result of R. Lazarsfeld \cite{laz} states the following:

\medskip
\noindent
{\it the induced morphism $f_*: H^i(\PP^n, \CC)\ra H^i(Y, \CC)$ is an
 isomorphism for  $i\le n-d+1$.}
\medskip 
 
 There is 
a natural vector bundle $\cE$ on $\PP^n$ of rank $d-1$ associated with
$f$, defined by the splitting 
$$f_*\cO_Y\cong \cO_{\PP^n}\oplus \cE^*$$
induced by 
the trace homomorphism. The map $f$ factors through
an embedding of $Y$ in the total space $|\cE |$ of $\cE$.
 Lazarsfeld proved that $\cE$ is always {\em ample}. 

This result implies in particular W. Barth's theorem \cite{barth}
for the cohomology of small codimensional subvarieties of projective
spaces.
Since Barth's fundamental paper, much attention has been 
paid to  
small codimension subvarieties of more general complex projective homogeneous
spaces \cite{sommese, faltings, bs}. This and the result of Lazarsfeld 
lead to the question 
of understanding  low degree branched coverings of such spaces $X=G/P$,
where $G$ is a semisimple complex Lie group and $P$ is a parabolic
subgroup. Note that the codimension
${\rm cod}\, (X, |\cE |)=d-1$.  
A recent paper of O. Debarre \cite{debarre} gives strong 
evidence to the following conjecture: 

\medskip\noindent {\bf Conjecture}. 
{\it Let $f: Y\ra X=G/P$ be a branched covering
where $Y$ is a nonsingular connected complex projective variety, 
and $Pic(X)=\ZZ$. Then the associated vector 
bundle $\cE$ on $X$ is ample.} 

\medskip The assumption $Pic(X)=\ZZ$ simply means that $P$ is a 
maximal parabolic subgroup. For $G=SL(n,\CC)$, the homogeneous space 
$G/P$ is then a grassmannian. The first half of the way to the
conjecture in that case was made in \cite{kim}, the second half in 
\cite{manivel}. The purpose of this paper is to prove the conjecture
for a few other cases. Our main result is the following:

\medskip\noindent {\bf Theorem A}. 
{\em Let $\LG_n$ be the lagrangian grassmannian
of maximal isotropic subspaces of a $2n$-dimensional symplectic vector
space. Then if $f: Y\ra \LG_n$ is a branched covering where
$Y$ is a nonsingular connected complex projective variety, 
then the associated vector bundle $\cE$ on $\LG_n$ is ample.}

\medskip Here, $\LG_n$ is homogeneous with  $G=Sp(2n,\CC),$
 the symplectic group;
with the exception of $\PP^{2n-1}$, $\LG_n$ is the only homogeneous space
for the symplectic group which is a hermitian symmetric space. 
Using Lazarsfeld's ideas, we immediately deduce from Theorem A the 
following Barth-Lefschetz type theorem:

\medskip\noindent {\bf Theorem B}. {\em Let $f: Y\ra \LG_n$ be a branched 
covering of degree $d$ where $Y$ is a nonsingular connected complex projective variety.
Then the induced morphism 
$$f^*: H^i(\LG_n,\CC)\ra H^i(Y,\CC)$$
is an isomorphism for $i\le {\rm dim}\,\LG_n-d+1$.} 

\medskip  The cohomological range of A. Sommese's Barth-Lefschetz type theorem  \cite{sommese2, sommese} concerning a
submanifold of a homogeneous space $X=G/P$ 
depends on  the $k$-ampleness (this is a suitable weakening of the
notion of ampleness) of the normal bundle of the
submanifold. 
To boot, by considering N. Goldstein's computation \cite{goldstein} of the $k$-ampleness 
of the tangent bundle of $X$, we can see that
 the codimension above
which one loses all cohomological information depends linearly on the 
rank of $G$, that is on the dimension of a maximal torus contained in
$G$; this
is reminiscent of G. Faltings' connectedness result for homogeneous spaces
(see \cite{faltings, goldstein} for the details). 
Here, our condition involves the dimension of $X=G/P$, which is in
general much larger, although in the case of $\PP^n$ this makes no difference. 

After the one of projective spaces, the easiest case should be 
the case of quadrics. Surprisingly enough, it does not seem to be so, 
despite partial results of the first author 
\cite{kim2}. Actually, we have been able to prove the conjecture 
only for quadrics of very small dimension:

\medskip\noindent {\bf Theorem C}. {\em Let $\QQ_n$ be the $n$-dimensional
nonsingular quadric with $3\le n\le 6$. 
If $f: Y\ra \QQ_n$ is a branched covering
where $Y$ is a nonsingular  connected complex projective variety, 
then the associated vector bundle $\cE$ on $\QQ_n$ is ample.}

\medskip Finally, let us mention that the conjecture by Debarre should be
extended to the cases of homogeneous spaces $X=G/P$ with $P$ not necessarily 
maximal. In this case, one cannot expect the vector bundle
associated with a branched covering to be ample. Nevertheless, it
will be $k$-ample for a suitable value $k$; this   
 implies Barth-Lefschetz type theorems. In section $5.3$ we propose 
a related conjecture  and give some evidence
to it by proving the following result:

\medskip\noindent {\bf Theorem D}. {\em Let $f: Y\ra X=G/P$ be a branched 
covering where $Y$ is a nonsingular connected complex projective variety, 
and $G=SL(n,\CC)$ or $Sp(2n,\CC)$. 
Then the associated vector bundle $\cE$ on $X$ is 
generated by its global sections.}

\medskip Actually, as in \cite{laz} we prove a stronger positivity property
for $\cE$ 
on $X=G/P$ with Picard
number one: $\cE\ot \cO_X(-1)$ is generated by its global
sections where $\cO_X(1)$ is the (very) ample generator of $Pic(X)$. 
To show this we set up an appropriate Castelnuovo-Mumford type criterion on $X$ and use vanishing theorems. 

In the case of lagrangian grassmannians, we have been obliged to get
out of the algebraic category, and instead we use vanishing theorems ``\`a
la Nakano'' for hermitian vector bundles. In section $2$, we
will recall a few elementary facts on the curvature of hermitian 
bundles, and state a vanishing theorem that will be suitable to our
purposes. 
With the help of a formula due to D.  Snow for the curvature of homogeneous
bundles, we will show in section $3$ how to apply this theorem 
 to the proof of Theorem A. In section $4$ we deal 
with low dimensional quadrics, where we make use of the so-called
spinor bundles to obtain Theorem C. 
Section $5$ is devoted to the proof of Theorem D, which is simple 
on ordinary flags, and much more involved in the symplectic case.

The first author thanks the Max-Planck-Institut f\"ur Mathematik in Bonn
for its warm hospitality during the preparation of part of this paper.

\section{On the curvature of hermitian vector bundles}
\subsection{Griffiths and Nakano positivity}

Let $E$ be a hermitian vector bundle on an $n$-dimensional 
nonsingular complex variety $X$.
The curvature $\Theta_E$ of the associated Chern connection is a 
$(1,1)$-form with values in the vector space of hermitian
endomorphisms of $E$. For every point $z\in X$ and every local
coordinate system $(z_j)_{1\le j\le n}$ at $z$,
the curvature can be written as 
$$\Theta_E=\sum_{jk\l\m}c_{jk\l\m}dz_j\wedge d\bar{z}_k\ot
e_{\l}^*\ot e_{\m},$$
with $c_{jk\l\m}=\bar{c}_{kj\m\l}$, where
$(e_1,\ldots ,e_r)$ is a local frame 
 of $E$ in a neighborhood of $z$. 
One may consider $\Theta_E$ as a hermitian form on $TX\ot E$:
$$\Theta_E(u)=\sum_{jk\l\m}c_{jk\l\m}(z)u_{j\l}\bar{u}_{k\m}$$
for $u=\sum_{j\l}u_{j\l}(\partial/{\partial z_j})\ot e_{\l}\in T_zX\ot E_z$. 

\medskip\noindent{\bf Definition}.  The hermitian vector bundle $E$
is Nakano (semi-)positive if $\Theta_E$ defines a (semi-)positive
definite hermitian form on $TX\ot E$. It is Griffiths
(semi-)positive if $\Theta_E$ is (semi-)positive
on non-zero decomposable tensors in $TX\ot E$. 

\medskip Nakano vanishing theorem states that the adjoint to a
 Nakano positive vector bundle on a compact K\"ahler
manifold has no higher cohomology:

\begin{theo}{\rm (cf. \cite{nakano})} 
Let $E$ be a hermitian vector bundle on a compact 
K\"ahler manifold. Suppose that $E$ is Nakano semi-positive, and that 
$\Theta_E$ is positive definite at least at one point. Then 
$$H^q(X,K_X\ot E)=0 \qquad \forall q>0.$$ \end{theo}

\subsection{A variant of Griffiths' vanishing theorem}

We will need in the sequel of this paper a consequence of the Nakano vanishing 
theorem, which can also be seen as a variant of Griffiths' vanishing
theorem \cite{griffiths}. We will state it in terms of {\em Schur powers}. 

Let $V$ be a finite dimensional complex vector space. To each
partition $\l$, i.e. to every non-increasing sequence of non-negative 
integers $\l_1\ge\cdots\ge \l_m\ge 0$ (=:the {\em parts} of $\l$), 
one associates a Schur
  power $S_{\l}V$, which is a  polynomial $Gl(V)$-module \cite{macdo}. 
This module reduces to zero if and only if the number of nonzero parts
of $\l$ (=:the {\em length} $l(\l)$) is larger than the dimension of
$V$. Otherwise, $S_{\l}V$ is an irreducible $Gl(V)$-module. 

The usual symmetric powers are a special case of  Schur powers.  They
correspond to partitions with a single nonzero part, or length one
partitions. Also the 
wedge powers correspond to partitions with all nonzero parts equal to
one. The tensor product of any Schur power with a symmetric or
a wedge power is described by the classical Pieri's rules
\cite{macdo}, which imply, in particular, that each Schur power 
$S_{\l}V$ is an irreducible component of the tensor product of 
$l(\l)$ symmetric powers (the exponents of which can be chosen as the 
parts of the partition $\l$). 

\begin{prop}\label{variant} Let $X$ be a compact K\"ahler manifold, 
and $E,F$ be 
hermitian vector bundles on $X$. Suppose that $E$ is Nakano 
semi-positive,  and that $\Theta_E$ is positive definite at least 
at one point. Suppose moreover that $F$ is Griffiths semi-positive. 
Then for any partition $\l$,
$$H^q(X,K_X\ot S_{\l}F\ot (\det F)^{l(\l)}\ot E)=0 
\qquad\forall q>0.$$ \end{prop}

\noindent {\em Proof.}
We first note that it is enough to treat the case of $l(\l )=1$,
that is, of symmetric powers. Indeed, if we can prove the vanishing
for one symmetric power, then by replacing $F$ by 
$F^{\oplus l}$ where $l=l(\l)$, which is also Griffiths semi-positive,
we deduce the vanishing for a tensor product of $l$ symmetric powers. Note
that  any
Schur power $S_{\l}E$ is a direct summand of such a tensor product.

Now let $\pi : Y={\mathbb P}(F^*)\lra X$ be the bundle of hyperplanes 
of $F$, and $\cO_F(1)$ be the tautological quotient line bundle on $Y$. 
By the isomorphisms \cite{griffiths}
$$H^q(X,K_X\ot S^lF\ot \det F\ot E)\cong H^q(Y,K_Y\ot \cO_F(l+r)\ot \pi^*E).$$
We want to show
that $\cO_F(m)\ot \pi^*E$ is Nakano semi-positive, and its curvature is positive
definite at least at one point, for any $m>0$, so that we can apply 
 Nakano vanishing theorem. Of course we endow this
vector bundle with two obvious metrics: on the one hand by the given 
metric on $E$ , which we simply pull back to one on $\pi^*E$, on 
the other hand by the pulled-back metric on $\pi^*F$, which provides
a metric on its quotient line bundle $\cO_F(1)$. Note that, 
since Griffiths semi-positivity is preserved by pull-backs and quotients, this
metric on $\cO_F(1)$ has Griffiths semi-positive curvature. Moreover, the 
Fubini-Study metric being positive on the projective space, it
has positive curvature on the fibers of $\pi$, 
that is, on ``vertical'' tangent vectors.  

Now the induced metric on $\cO_F(m)\ot \pi^*E$ has curvature
$$\Theta_{\cO_F(m)\ot \pi^*E} =
m\Theta_{\cO_F(1)}\ot Id_{\pi^*E}+\pi^*\Theta_E.$$
Let us choose a local frame $(e_1,...,e_r)$ of $E$ in a neighborhood of
$x=\pi (y)$ where $y\in Y$. Then if 
$u=\sum_{i=1}^ru_i\ot e_i\in T_yY\ot E_{\pi(y)}$, we have  
$$\Theta_{\cO_F(m)\ot \pi^*E}(u)=
m\sum_{i=1}^r\Theta_{\cO_F(1)}(u_i)+\Theta_E(\pi_*u).$$

This is a sum of non-negative terms. Moreover, suppose that $\Theta_E$
is positive definite at $x=\pi(y)$. Then if $\Theta_{\cO_F(m)\ot \pi^*E}(u)$
is equal to zero, we must have $\pi_*u=0$. 
This means that $u$  has only vertical components, and since 
$\Theta_{\cO_F(1)}$ is positive definite on the fibers of $\pi$, we 
must have $u_i=0$ for each $i$, hence $u=0$. Thus $\Theta_{\cO_F(m)\ot
\pi^*E}$ is positive definite at $y$, which proves our claim. \qed

\medskip 
We will use  the following vanishing theorem in the sequel, which is a slight
generalization of Le Potier vanishing theorem \cite{lepotier}:

\begin{prop}\label{lp} Let $E_1,\ldots ,E_m$ be globally generated vector
  bundles of ranks $e_1,\ldots ,e_m$, and $L$ be an ample line bundle, 
on a nonsingular  complex projective variety $X$. Then
$$H^q(X,K_X\ot\we^{e_1-k_1}E_1\otc \we^{e_m-k_m}E_m\ot L)=0$$
for $q>k_1+\cdots +k_m$. \end{prop}

\subsection{The curvature of homogeneous bundles}

We now consider the very specific case of a homogeneous vector bundle 
$E=G\times_P E_0$ of rank $r$ 
on a generalized flag manifold $X=G/P$, where $P$ is
a parabolic subgroup of a semisimple complex Lie group $G$. Here $E_0$
is a finite dimensional $P$-module. 

Assume that $E$ is generated by its global sections
so that the evaluation $P$-module homomorphism
$$\phi : V=H^0(X, E)\ra E_0$$ 
is surjective. For a given 
hermitian metric on $V$,  we have an induced metric on the fiber $E_0$,
and a left-invariant hermitian metric on $E$. 

The curvature of this metric was computed by D. Snow \cite{snow}. We need
some notation to state his formula. Firstly, let us fix a maximal torus
and a Borel subgroup in $G$, from which we deduce a root system for
its Lie algebra. Up to conjugacy, the parabolic subgroup $P$ of $G$ is
then generated by  this Borel subgroup, and a closed set of 
 certain positive roots. The opposite of these 
roots form a set $I$ of simple roots, and $\Phi_X$, called the roots of $X$,
 is then defined as the 
set of positive roots which, when decomposed as a linear combination
of simple roots, have positive coefficients on the roots in $I$. 

Let  $v_1, \ldots , v_m$ be a basis
of root vectors for $V$ and $e_1,\ldots ,e_r$ be a basis for $E_0$, such that
$$\phi (v_k)=\left\{ \begin{array}{ll} e_k & {\rm for}\; 1\le k\le r,
    \\ 0 & {\rm for}\; r<k\le m. \end{array}\right.$$

\begin{prop}\label{snow} 
The curvature of the hermitian metric induced on $E$ by
    the evaluation morphism $\phi$, is given at the
identity coset by the following formula:
$$\Theta_E =\sum_{\alpha,\beta\in\Phi_X}\sum_{k>r}
dx_\alpha\we d\bar{x}_\beta\ot
\phi(X_\alpha \cdot v_k)\ot\overline{\phi(X_\beta\cdot v_k)^t},$$ 
where transposition is taken with respect to the given basis. 
\end{prop}

\noindent
{\em Remark.}
Being induced by a constant metric on a trivial vector bundle, 
this metric on $E$ is Griffiths semi-positive.
But in general, is neither Griffiths positive, nor Nakano
semi-positive and a fortiori not even Nakano positive. 
It would be useful to find conditions on $E_0$ ensuring that 
the induced metric on $E$ has some of these properties. 

\section{Symplectic grassmannians}

In this section, we consider a symplectic grassmannian $\LG_n$,
defined as the subvariety of the usual grassmannian $\GG_{n,2n}$,
consisting of maximal isotropic spaces in $\CC^{2n}$ 
with respect to a fixed symplectic form. For example, $\LG_1=\PP^1$, and 
$\LG_2$
is a hyperplane section of $\GG_{2,4}$ in its Pl\"ucker embedding, hence a
three dimensional quadric. In general, $\LG_n$ is a hermitian symmetric
space of dimension $n(n+1)/2$, homogeneous under the action of the
symplectic group $Sp(2n,\CC)$. If we write $\LG_n=Sp(2n,\CC)/P$, then
the parabolic subgroup $P$ is defined (up to conjugacy) by the 
(unique) long simple root.

\subsection{A Castelnuovo-Mumford criterion}

We will need a Castelnuovo-Mumford criterion on the symplectic 
grassmannian $\LG_n$, that is, a cohomological criterion for a 
coherent sheaf on $\LG_n$ to be generated by its global sections. 
On the usual grassmannian, such a criterion is deduced in \cite{kim} and
\cite{manivel}
from the simple fact that  a general section of 
a some homogeneous bundle (precisely, the sum
of a suitable number of copies of the tautological quotient vector bundle)
vanishes precisely at one point. Unfortunately, we cannot give the same 
argument on symplectic grassmannians (see \cite{fulton} for a discussion 
of closely related questions). Moreover, previously mentioned 
Castelnuovo-Mumford criterion on $\GG_{n,2n}$, which we may also use 
on the subvariety $\LG_n$, will not be good enough for our purposes. 

Let $Q_n$ be the tautological quotient bundle of rank $n$ 
on $\LG_n$. We define 
$$\cL_i^{n,j}=\left\{ \begin{array}{lll} \we^iQ_n^* & {\rm if} & i<j, \\
\cO_{\LG_n}(-1) & {\rm if} & i=j, \\ 0 & {\rm if} & i>j. \end{array}\right.$$
Our Castelnuovo-Mumford criterion on $\LG_n$ will be the following:

\begin{prop} Let $\cF$ be a coherent sheaf on $\LG_n$. Then $\cF$ 
is generated by its global sections as soon as 
$$H^q(\LG_n,\cL^{n,1}_{i_1}\otc \cL^{n,n}_{i_n}\ot\cF)=0$$
for $q\ge i_1\pc i_n>0$ and $q>i_1\pc i_n=0$. \end{prop}

\noindent {\em Remark}. In the following proof and in the
sequel, we will repeatedly use the following simple 
observation: for a coherent sheaf $\cG$  on a variety $X$, 
let $\cR_{\bullet}\ra\cG\ra 0$ be a finite resolution of $\cG$ indexed by 
non-negative integers. Then 
the cohomology group $H^q(X,\cG)$ vanishes as soon as 
$$H^{q+m}(X,\cR_m)=0$$
for all $m\ge 0$. This can be  easily seen by breaking the resolution
into short exact sequences, and inspecting the corresponding long
exact sequences of cohomology groups. Similarly, if we have  a resolution 
$0\ra\cG\ra\cS_{\bullet}$, then  $H^q(X,\cG)$ vanishes as soon as 
$H^{q-m}(X,\cS_m)=0$ for all $m\ge 0$.
\bigskip

\noindent {\em Proof.} We proceed by induction on $n$. The 
space of global sections of $Q_n$ is $\CC^{2n}$, and the zero-locus
of the section of $Q_n$ corresponding to some non-zero vector $v$ is the 
space of maximal isotropic subspaces of $\CC^{2n}$ containing $v$. 
This locus can be identified with the lagrangian grassmannian
$\LG_{n-1}$ of maximal isotropic spaces for the induced symplectic form
on $v^{\perp}/v\simeq\CC^{2n-2}$. Moreover, these loci 
cover $\LG_n$ $v$ as varies. To prove that $\cF$ is generated by its global
sections on $\LG_n$, it is therefore enough to check that its
restriction to any such $\LG_{n-1}$ is generated by its global sections, 
and the surjectivity of the restriction map
$$H^0 (\LG_n, \cF) \lra H^0 (\LG_{n-1}, \cF|_{\LG_{n-1}}).$$

By the induction hypothesis, $\cF|_{\LG_{n-1}}$ will be generated by its global
sections as soon as 
$$H^q(\LG_{n-1},\cL^{n-1,1}_{i_1}\otc \cL^{n-1,n-1}_{i_{n-1}}
\ot\cF|_{\LG_{n-1}})=0$$
for $q\ge i_1\pc i_{n-1}>0$ and $q>i_1\pc i_{n-1}=0$. 
From the short exact sequence 
$$0\ra \cO_{\LG_{n-1}} \ra Q_n^*|_{\LG_{n-1}} \ra Q_{n-1}^* \ra 0,$$
we get long exact sequences
$$0\ra \cO_{\LG_{n-1}} \ra Q_n^*|_{\LG_{n-1}} \ra
\wedge^2 Q_n^*|_{\LG_{n-1}} \ra \cdots \ra 
\wedge^kQ_n^*|_{\LG_{n-1}} \ra \wedge^kQ_{n-1}^* \ra 0,$$
from which we deduce that the proceeding vanishing will hold if
$$H^q(\LG_{n-1}, (\cL^{n,1}_{j_1}\otc \cL^{n,n-1}_{j_{n-1}}
\ot\cF)|_{\LG_{n-1}})=0$$
for $q\ge j_1\pc j_{n-1}>0$ and $q>j_1\pc j_{n-1}=0$. 
Now we  compare cohomology groups of coherent sheaves on $\LG_n$, and
of their restriction to $\LG_{n-1}$, thanks to the exact Koszul complex of the
section of $Q_n$ defining $\LG_{n-1}$: 
$$0\ra \cO_{\LG_n}(-1) \ra \we^{n-1}Q_n^* \ra \cdots \ra 
Q_n^* \ra \cO_{\LG_n}\ra \cO_{\LG_{n-1}}\ra 0.$$
The above vanishing on $\LG_{n-1}$ will thus be a consequence of 
the following ones on $\LG_n$:
$$H^q(\LG_n,\cL^{n,1}_{j_1}\otc \cL^{n,n-1}_{j_{n-1}}\ot \cL^{n,n}_{j_n}
\ot\cF)=0$$
for $q\ge j_1\pc j_n>0$ and $q>j_1\pc j_n=0$. This is precisely  
our criterion. Finally, again, by the use of the above Koszul complex, 
the surjectivity of the restriction map 
$$H^0 (\LG_n, \cF) \lra H^0 (\LG_{n-1}, \cF|_{\LG_{n-1}})$$
is a consequence of the vanishing of $H^q(\LG_n,\cL^{n,n}_q\ot\cF)$
for $q>0$, which is also included in our criterion. \qed

\subsection{The key observation}

If we try to apply the previous Castelnuovo-Mumford criterion on $\LG_n$ to
$\cE(-1)$ where $\cE$ is the vector bundle associated to some finite
covering, the classical vanishing theorem of Griffiths or Le Potier
will not suffice. Actually, we will need our variant of Griffiths 
vanishing \ref{variant}, which involves  a Nakano positive vector
bundle. Our key observation is the following: 

\begin{prop} The vector bundle $Q_n(1)$ on $\LG_n$ is Nakano 
positive. \end{prop} 

\noindent{\em Remark}. Note that the corresponding statement 
for the tautological quotient bundle $Q$ on a usual grassmannian is false. 
However, since $Q$ is a quotient bundle of a trivial bundle which is
Griffiths semi-positive, it follows from general properties 
of hermitian bundles that $Q(1)=Q\ot\det Q$ is Nakano semi-positive. 

\bigskip\noindent {\em Proof.} Let $ J=\left( \begin{array}{cc} 0 & I_n\\
-I_n & 0 \end{array} \right) $ be the $2n\times 2n$ matrix associated with
a symplectic form on $\CC^{2n}$,
where $I_n$ is the identity matrix. We follow the conventions and
notations of \cite{bourbaki}. In particular, we denote by $(e_1,\ldots
,e_n,e_{-n},\ldots ,e_{-1})$ the canonical basis of $\CC^{2n}$.
We have the roots of $\LG_n$,
$$\Phi_{\LG_n}=\{\e_i+\e_j:\;1\le i\le j\le n\},$$
and the the following corresponding root vectors:
\begin{eqnarray*}
X_{ij}=X_{\e_i+\e_j}=E_{i,-j}+E_{j,-i} & \mbox{\, for\, } 1\le i<j\le
n, \\ X_{ii}=X_{2\e_i}=E_{i,-i} & \mbox{\, for\, } 1\le i\le n.
\end{eqnarray*}
Let us apply Snow's formula \ref{snow}
to $Q_n$, at the point of $\LG_n$ given by 
the maximal isotropic subspace of $\CC^{2n}$ generated by $e_1,\ldots
,e_n$.  The vectors $v_m=e_{-m}$, for $1\le m\le n$, form a basis of 
the kernel of the evaluation map, and we have
\begin{eqnarray*}\sum_{1\le i\le j\le n} \phi (X_{ij}\cdot v_m)
\, dx_{ij} & = & \sum_{1\le i<j\le n}(\delta_{jm}e_i+\delta_{im}e_j)
\, dx_{ij}+e_m\, dx_{mm} \\ & = & \sum_{i=1}^n e_i\, dx_{im}\end{eqnarray*}
by adopting the notation $dx_{ij}=dx_{ji}.$
We then obtain the curvature forms
\begin{eqnarray*}
\Theta_{Q_n} & = &
\sum_{1\le i,j,k\le n}  dx_{ik}\,\we\, d\bar x_{jk}\ot e_i\ot e_j^*,\\
\Theta_{\cO_{\LG_n}(1)} & =  & {\rm Trace}\, \Theta_{Q_n} 
= \sum_{1\le i,k\le n} dx_{ik}\,\we\, d\bar x_{ik}.
\end{eqnarray*}
The second of these formulae is a consequence of the first 
since $\cO_{\LG_n}(1)=\det Q_n$. Then we have the hermitian form
$$\Theta_{Q_n(1)} (u)  = 
\sum_{1\le i,j,k\le n} (\, u_{jk}^i \,\overline{u}_{ik}^j + |u_{ik}^j|^2\, )$$
where $u_{jk}^i=u_{kj}^i$. We can rewrite this form as 
$$\Theta_{Q_n(1)}(u)=\sum_{i<j<k}\Theta_{ijk}(u)
+\sum_{i\ne j}\Theta_{ij}(u)+\sum_i\Theta_i(u),$$
where the different terms of this sum are hermitian forms 
involving disjoint sets of components of $u$. Specifically,
$$\begin{array}{rcl}
\Theta_{ijk}(u) & = & 2\,(|u_{jk}^i|^2+|u_{ik}^j|^2+|u_{ij}^k|^2)+ \\ 
& & +u_{jk}^i\bar{u}_{ik}^j + \bar{u}_{jk}^i u_{ik}^j
+u_{jk}^i\bar{u}_{ij}^k+\bar{u}_{jk}^i u_{ij}^k+
u_{ik}^j\bar{u}_{ij}^k +\bar{u}_{ik}^j u_{ij}^k \\
 & = & |u_{jk}^i+u_{ik}^j+u_{ij}^k|^2+|u_{jk}^i|^2+|u_{ik}^j|^2+
|u_{ij}^k|^2. 
\end{array}$$
It is a positive definite hermitian form in the three variables 
$u_{ij}^k, u_{jk}^i$ and $u_{ik}^j$, while
$$\Theta_{ij}(u)\,= \,|u_{jj}^i+u_{ij}^j|^2+2\, |u_{ij}^j|^2, \quad
\Theta_i(u)\, = \,2\, |u_{ii}^i|^2$$
are positive definite hermitian forms in $u_{ij}^j,u_{jj}^i$
and $u_{ii}^i$, respectively. 
This clearly implies that the hermitian form $\Theta_{Q_n(1)}$
is positive definite, and the proposition is proved. \qed

\subsection{A lemma on tensor products}

Before proving Theorem A, we need the following simple remark on 
Schur powers, which is an easy consequence of Pieri's rules:

\begin{lemm}\label{schur} 
Let $S_{\l}V$ be an irreducible component of a tensor 
product of $m$ symmetric powers $S^{k_1}V\otc S^{k_m}V$. Then 
the length of $\l$ is at most $m$, and if it is equal to $m$, 
$S_{\l}V$ is also an irreducible component of the tensor 
product $\we^mV\ot S^{k_1-1}V\otc S^{k_m-1}V$.\end{lemm}

\noindent {\em Remark.} More precisely, $S_{\l}V$ has the
same multiplicity in both tensor products, but we do not need this fact in this
context. 
  
\bigskip\noindent {\em Proof.}
Because of Pieri's rules, $S_{\l}V$ is a
component of $S^{k_1}V\otc S^{k_m}V$ if and only if there exists
a semistandard tableau $T$ of shape $\l$ and weight
$\kappa=(k_1,...,k_m)$: that is, a numbering of the diagram of 
the partition $\l$, weakly increasing on rows, strictly increasing on
columns such that each integer $i$ occurs $k_i$ times; see \cite{macdo}. 

The length of $\l$ is the length of its first column. Since this
column is numbered in a strictly increasing way by integers not 
exceeding $m$, its length is certainly bounded by $m$. If it is 
equal to $m$, this column must be numbered by $1,\ldots ,m$. If we remove
the first column, we get a new semistandard tableau $S$, of shape $\mu$ 
and weight $(k_1-1,\ldots ,k_m-1)$. Since $S_{\l}V$ is an irreducible 
component of $\we^mV\ot S_{\mu}V$, the lemma is proved. \qed 

\subsection{Proof of Theorem A}

Let $f : Y\lra \LG_n$ be a degree $d$ covering, and $\cE$ be the associated 
vector bundle on $\LG_n$ of rank 
$d-1$.  We want to show that 
$\cE(-1)$ is  generated by its global sections. Using Serre duality
and the fact that $K_{\LG_n}=\cO_{\LG_n}(-n-1)$, we see that 
our Castelnuovo-Mumford criterion is equivalent to the vanishings of 
$$H^{N-q}(\LG_n,(\cL^{n,1}_{i_1}\otc \cL^{n,n}_{i_n})^*
\ot\cE^*(-n))$$
for $q\ge i_1\pc i_n>0$ or $q>i_1\pc i_n=0$, where $N=\dim\,\LG_n$.
Because of the definition of $\cE$, Serre duality again and Leray spectral
sequence, 
these vanishings are equivalent to the isomorphisms  
$$H^q(\LG_n,K_{\LG_n}\ot\cL^{n,1}_{i_1}\otc \cL^{n,n}_{i_n}(n))
\cong H^q(Y,K_{Y}\ot f^*(\cL^{n,1}_{i_1}\otc \cL^{n,n}_{i_n}(n)))$$
for $q\ge i_1\pc i_n>0$ or $q>i_1\pc i_n=0$. 

If $i_k=k$ for $1\le k\le n$, and $q=i_1\pc i_n=N$,
both cohomology groups are equal to $\CC$. If $i_k=0$ for $1\le k\le n$, then
both groups vanish by  Kodaira vanishing theorem. We shall 
show that all the other groups also vanish. 

We first deal with cohomology groups on $\LG_n$. Note that they 
could be computed using representation theory and applying 
the Borel-Weil-Bott theorem, but we will give a quite different 
argument mainly using  Griffiths' vanishing theorem and its variant
in section
2.2. 
 
Let us consider a given $n$-tuple $i_1,\ldots ,i_n$, with 
$0\le i_k\le k$. Let $a$ be the number of indices $k$ such
that $i_k<k$, which is supposed to be positive. Note that 
the sum of the other $i_k$'s is  bigger than or
equal to $(n-a)(n-a+1)/2$. 
It will therefore be enough to prove that 
$$H^q(\LG_n,K_{\LG_n}\ot\we^{l_1}Q_n^*\otc \we^{l_a}Q_n^*(a))=0$$
for $q\ge l_1\pc l_a+(n-a)(n-a+1)/2$, and for all sequences 
$l_1,\ldots ,l_a$ with $l_k<n-a+k$. (This condition is equivalent 
to the existence of a strictly increasing sequence $1\le m_1<\cdots
<m_a\le n$ such that $l_k<m_k$ for all $k$.) 

We then use the tautological exact sequence $0\ra Q_n^*\ra V\ra Q_n\ra
0$ on $\LG_n$ where $V=\CC^{2n}$ (identified with the corresponding trivial
bundle on $\LG_n$; recall that if $W$ is maximal isotropic subspace of $V$, 
then the quotient $V/W$ is canonically identified with the dual
$W^*$), and its induced long exact sequences of wedge powers 
$$0\ra \we ^l Q_n^*\ra \we^lV\ra \cdots \ra \we ^{l-m}V\ot S^mQ_n \ra 
\cdots\ra S^l Q_n\ra 0.$$

{\bf First case: $a<n$}. Using the above long exact sequences for each 
wedge power, we are reduced to prove that 
$$H^{q-m_1-\cdots -m_a}(\LG_n,K_{\LG_n}\ot S^{m_1}Q_n\otc S^{m_a}Q_n(a))=0$$
for $q\ge l_1\pc l_a+(n-a)(n-a+1)/2$, and for all sequences 
$m_1,\ldots ,m_a$ with $m_k\le l_k$ for all $k$. 

Let $S_{\l}Q_n$ be an irreducible component of $S^{m_1}Q_n\otc S^{m_a}Q_n$. 
If $l(\l)< a$, then the vanishing follows from
Griffiths' vanishing theorem, more precisely from its variant
\ref{variant} by setting $F=Q_n$, which is Griffiths semi-positive 
as a quotient of a trivial bundle, and by setting 
$E=\cO_{\LG_n}(1)$, which is positive. If $l(\l )=a$, then 
lemma \ref{schur} implies that it is enough to prove that 
$$H^{q-m_1-\cdots -m_a}(\LG_n,K_{\LG_n}\ot S^{m_1-1}Q_n\otc S^{m_a-1}Q_n
\ot\we^aQ_n(a))=0.$$
Because $\we^aQ_n=\we^{n-a}Q_n^*(1)$,  we can use the above complexes
again to reduce to the vanishings of 
$$H^{q-m_1-\cdots -m_a-m_{a+1}}(\LG_n,K_{\LG_n}\ot S^{m_1-1}Q_n\otc S^{m_a-1}Q_n
\ot S^{m_{a+1}}Q_n(a+1)),$$
with $m_{a+1}\le n-a$. We then use the same argument,
repeated $n-a-1$ times, and we are finally reduced to prove the vanishings of 
$$\begin{array}{l}
H^{q-m_1-\cdots-m_{n-1}}(\LG_n,K_{\LG_n}\ot S^{m_1-n+a+1}Q_n\otc 
S^{m_a-n+a+1}Q_n\ot \\
\hspace*{5cm}\ot S^{m_{a+1}-n+a+2}Q_n\otc S^{m_{n-1}}Q_n(n-1)),
\end{array}$$
where $m_{a+k}\le n-a-k+1$. 
Firstly note that $q-m_1-\cdots-m_{n-1}$ is positive since 
$m_1\pc m_a\le l_1\pc l_a$ and $m_{a+1}\pc m_{n-1}<(n-a)(n-a+1)/2$.
Recall moreover that $l_1\le n-a$, so that $m_1-n+a+1\le 1$. 
If it is equal to zero, we can apply Griffiths' vanishing 
theorem as we did above. 
If it is equal to one, we can factor out a $Q_n(1)$ and apply 
\ref{variant} by setting $F=Q_n$ and $E=Q_n(1)$.

\smallskip {\bf Second case: $a=n$}. We need to prove that 
$$H^q(\LG_n,K_{\LG_n}\ot\we^{l_1}Q_n^*\otc \we^{l_n}Q_n^*(n))=0$$
for $q\ge l_1\pc l_n>0$, and for all sequences 
$l_1,\ldots ,l_n$ with $l_k<k$. Note that in particular 
$l_1=0$. Let us choose $p$ such that $l_p>0$. 
We then use our long exact sequences for each 
wedge power other than $\we^{l_p}Q_n^*$. We are thus reduced to prove that 
$$H^q(\LG_n,K_{\LG_n}\ot S^{m_2}Q_n\otc S^{m_n}Q_n\ot\we^{n-l_p}Q_n(n-1))=0$$
for $q\ge l_p$, and for all sequences 
$m_1,\ldots ,m_{p-1},m_{p+1},\ldots,m_{n}$ with $m_k\le l_k$.
 Note that we have only $n-2$ symmetric powers in this 
tensor product. Let $S_{\l}Q_n$ be some component of it. If $l(\l)< n-2$, 
or $l(\l)=n$, we can apply Griffiths'
vanishing theorem. If $l(\l)=n-1$, than the last non-zero part of the
partition $\l$ must be equal to one. Let $\mu$ be the partition
obtained by deletion 
of this last part: then $l(\mu)=n-2$, and $S_{\l}Q_n$ is 
a component of $S_{\mu}Q_n\ot Q_n$. We can therefore apply 
\ref{variant} by setting $F=Q_n$ and $E=Q_n(1)$ to get the required vanishing. 

\smallskip There remains to prove the similar vanishings on the 
finite covering $Y$ of $\LG_n$. But the very same argument ``pulls
back'' to $Y$. Indeed, we used nothing more on $\LG_n$ than the 
tautological exact sequence, the Griffiths semi-positivity of $Q_n$
and the Nakano positivity of $\cO_{\LG_n}(1)$ and $Q_n(1)$. 
But the tautological sequence is clearly preserved on $Y$, as well as the 
Griffiths and Nakano semi-positivity of pulled-back bundles. Strict
positivity, of course, is not preserved by pull-back. But for a finite
covering, it is preserved outside the ramification locus, since
outside this locus we have a local isomorphism for the complex
topology.  Since only positivity at one point is required
 to apply the vanishing theorem
\ref{variant}, we conclude our theorem. \qed
\medskip

 By proceeding  as in   
\cite{laz2}, we obtain the following homotopy result with a weaker homotopy
bound than the cohomology bound in Theorem B.

\begin{prop}
 Let $Y$ and $f$ be as in Theorem A. Then for any point $y\in Y$,
$$f_*:\pi_i (Y,y)\ra \pi_i (\LG_n, f(y))$$
 is an isomorphism for
$i\le \dim \LG_n - \, {\rm max}\, \{ d-1, k+1\}$ where $k$ is the 
$k$-ampleness of the tangent bundle $T_{\LG_n}$.\end{prop}

\medskip  
\noindent {\em Proof}. Let $H:=\LG_n$.
 Consider $f\times f:Y\times Y \ra H\times H$, 
and let $\Delta_Y$ and $\Delta_H$ be the diagrams in 
$Y\times Y$ and $ H\times H$,
repectively.
Using the
fact that $\cE (-1)$  is generated by its global sections on $H$,
we see that 
the inclusion $ \Delta_Y\hookrightarrow X:=(f\times f)^{-1}(\Delta_H)$
 induces  a morphism 
$$H_i (Y,\ZZ)\ra H_i(X,\ZZ)$$
which is an isomorphism for  
$i\le \dim H -d$  and is surjective for $i=\dim H-d+1$ 
 (cf. \cite[Lemma 5.4]{laz2}).
Note that $Y$ is simply connected if $ \dim H \le d$ 
\cite[Corollaire 2.6]{debarre}.
 Now we consider the following commutative diagram of
$\ZZ$-homology
 sequences 
induced by $f\times f$
$$\begin{array}{cccccccccc}
\ra\!&\!H_{i+1}(Y\times Y,X)\!&\!\ra\!&\!H_i (X) \!&\!\ra  \!&\!H_i (Y\times Y) \!&\!\ra \!&\! H_i (Y\times Y, X)\!&\!\ra \\
&\downarrow&&\downarrow & & \downarrow & &\downarrow&&\\
\ra\!&\!H_{i+1}(H\times H,H)\!&\!\ra\!&\!H_i (\Delta_H) \!&\!\ra\!&\!H_i (H\times H)\!&\!\ra\!&\! H_i (H\times H,\Delta_H) \!&\!\ra 
\end{array}$$ 
The first vertical homomorphism is an isomorphism 
if $i+1\le \dim H -k$ by a version of Barth-Lefschetz theorem
due to Sommese-Van de Ven \cite{sv}. 
 Since $H$ has homology groups
in even degrees only
and it doesn't have torsion,   the induced morphism
$$f_*: H_i(Y,\ZZ)\ra H_i (H,\ZZ)$$
is an isomorphism  for 
$i\le \dim H -\, {\rm max}\,\{d-1, k+1\}.$ By J.H.C.
Whitehead's theorem \cite{whitehead}, we conclude the proposition.
\qed

\section{Quadrics}

\subsection{Spinor bundles}

The $n$-dimensional nonsingular quadric $\QQ_n$  is a 
homogeneous space $SO(n+2,\CC)/P$, where $P$ is the maximal 
parabolic subgroup associated, in the notation of \cite{bourbaki},
with the first simple root of $so(n+2,\CC)$. It is a simple 
Lie algebra of type $B$ if $n$ is odd,  and of type $D$ if  $n$ is even. 
Note that the semisimple part of $Lie(P)$ is isomorphic to 
$so(n,\CC)$. If $n=2m+1$, this semisimple Lie algebra has a spinor 
representation which defines a homogeneous bundle $S$ on $\QQ_n$
of rank $2^m$. If $n=2m$, there are two non isomorphic
half spinor representations, giving rise to two homogeneous vector
bundles $S'$ and $S''$ of rank $2^{m-1}$. 

These {\em spinor bundles} have been extensively studied by 
Ottaviani \cite{ottaviani}. 
Among other things, he showed the following:

\noindent 
On $\QQ_{2m+1}$, we have a tautological exact sequence
$$0\ra S\ra V\ra S(1)\ra 0$$
and an isomorphism $S^*\cong S(1)$
where $V$ is trivial of rank $2^{m+1}$.

\noindent 
On $\QQ_{2m}$, we have
two tautological exact sequences (dual to each other)
$$0\ra S'\ra U\ra S''(1)\ra 0,\qquad 
0\ra S''\ra U^*\ra S'(1)\ra 0,$$
where $U$ is trivial of rank $2^m$. 
If $m\equiv 0 \; ({\rm mod}\; 2)$ then $S'^*\cong S'(1)$ and
$S''^*\cong S''(1)$, and if $m\equiv 1 \; ({\rm mod}\; 2)$ then
$S'^*\cong S''(1)$ and $S''^*\cong S' (1)$.

\noindent
Furthermore, spinor bundles behave 
nicely under restriction to generic hyperplane sections: namely, the restriction
of $S$ to $\QQ_{2m}$ is $S'\op S''$, while the restriction of $S'$ or
$S''$ to $\QQ_{2m-1}$ is the spinor bundle on that quadric.

\subsection{Quadrics of dimension five}

On $\QQ_5$, we have the following Castelnuovo-Mumford criterion: 

\begin{prop}\label{cmq5} 
A coherent sheaf $\cF$ on $\QQ_5$ is generated by its global sections as soon as 
$$\begin{array}{lcll}
H^{i+j}(\QQ_5,\cF(-l))=0 & {\rm for} & (i,l)=(1,1),(2,2),(5,3), & j\ge 0,\\
H^{i+j}(\QQ_5,\cF\ot S(-l))=0 & {\rm for} & (i,l)=(1,0),  & j\ge 0,
\quad {\rm and}\\
H^i(\QQ_5,\cF\ot \we^2S(-l))=0 & {\rm for} & (i,l)=(3,1),(4,2). & 
\end{array}$$ \end{prop}

\noindent {\em Proof}. Let $x$ be a point in $\QQ_5$, and $\QQ_4$ a 
four dimensional quadric in $\QQ_5$ containing $x$. Since $\QQ_4$ is a
hyperplane section, we have an exact sequence 
$$0\ra\cO_{\QQ_5}(-1)\ra\cO_{\QQ_5}\ra\cO_{\QQ_4}\ra 0,$$
showing in particular that the restriction map 
$$H^0(\QQ_5,\cF)\lra H^0(\QQ_4,\cF|_{\QQ_4})$$
is surjective, since $H^1(\QQ_5,\cF(-1))=0$ by hypothesis. 

On $\QQ_4$, we have two half spinor bundles $S'$ and $S''$ of rank 2, 
and $x$ can be realized as the zero-locus of a section of $(S'\op S')^*$. 
Indeed, we can identify $\QQ_4$ with the grassmannian $\GG_{2,4}$
in such a way that $S'^*$ identifies with the tautological quotient 
bundle on that grassmannian. 
Using the Koszul complex  associated with a section of $(S'\op S')^*$, 
we see that $\cF|_{\QQ_4}$ is globally
generated at $x$ as soon as 
$$H^q(\QQ_4,\we^q(S'\op S')\ot\cF|_{\QQ_4})=0\qquad \forall q>0.$$
Now recall that $S|_{\QQ_4}=S'\op S''$, so that for $q=1$, this
vanishing is a consequence of the following conditions on $\QQ_5$:
$$H^1(\QQ_5,\cF\ot S)\cong H^2(\QQ_5,\cF\ot S(-1))=0,$$
where the isomorphism is obtained by the tautological exact sequence on
$\QQ_5$ and our hypothesis.
For $q=2$, we want to show that $H^2(\QQ_4,S'\ot S'\ot\cF|_{\QQ_4})=0$.
Tensoring the second tautological sequence
by $S'(-1)\ot\cF|_{\QQ_4}$ on $\QQ_4$, we reduce to the vanishings of 
$H^2(\QQ_4,S'(-1)\ot\cF|_{\QQ_4})$ and $H^3(\QQ_4,S'\ot
S''(-1)\ot\cF|_{\QQ_4})$. The first of these groups is dealt with as
above. For the second one, we note that $S'\ot S''$ is a summand of 
$\we^2S|_{\QQ_4}$, so that our vanishing condition is a consequence of
$$H^3(\QQ_5,\cF\ot \we^2S(-1))=H^4(\QQ_5,\cF\ot \we^2S(-2))=0.$$
Finally, the cases $q=3$ and $q=4$ are similar to that of $q=1$. \qed 

\medskip\noindent {\em Proof of Theorem C for $3\le n\le 5$}.
Recall that $\QQ_3\simeq\LG_2$ and $\QQ_4\simeq\GG_{2,4}$,
so that Theorem C follows for $n=3$ from Theorem A,
and for $n=4$ from \cite{manivel}. Let now $n=5$. 

Let $f: Y\ra\QQ_5$ be a branched covering and $\cE$ be the associated vector
bundle. We prove that $\cE(-1)$ is generated by global sections, using
our Castelnuovo-Mumford criterion on $\QQ_5$. Using Serre duality,
the definition of $\cE$ and the fact that
$K_{\QQ_5}=\cO_{\QQ_5}(-5)$, we see that what we need to show is 
that certain cohomology groups of type
$$H^i(\QQ_5,K_{\QQ_5}\ot \we^jS(4-l))$$
are equal to the corresponding cohomology groups of the pulled-back
bundles on $Y$. We show that all these groups vanish, with one
exception. 

This is clear for $j=0$ since $l$ is at most two, which allows to use 
Kodaira vanishing. For $j=1$, we can use the tautological exact sequence 
$0\ra S\ra V\ra S^*=S(1)\ra 0$, which easily implies that 
$$H^i(\QQ_5,K_{\QQ_5}\ot S(m))\cong H^{i-n}(\QQ_5,K_{\QQ_5}\ot S(m+n)),$$
provided that $m,m+n$ and $i-n$ are positive. Since $S=S^*(-1)$, 
$\det S^*=\cO_{\QQ_5}(2)$ and $S^*$ is globally generated, the 
Griffiths vanishing theorem then gives 
$$H^i(\QQ_5,K_{\QQ_5}\ot S(4-l))=
H^1(\QQ_5,K_{\QQ_5}\ot S^*\ot \det S^*(i-l))=0$$
for $0\le l<i\le 4$, which is the case in our criterion. 

Finally for $j=2$, we have two conditions. The first one concerns 
$$H^3(\QQ_5,K_{\QQ_5}\ot \we^2S(3))\cong H^3(\QQ_5,K_{\QQ_5}\ot \we^2S^*(1))$$
which is zero by the Le Potier vanishing theorem. For the second one, 
we shall prove that 
$$H^4(\QQ_5,K_{\QQ_5}\ot \we^2S(2))\cong\CC.$$
We first notice that $S^*$ is globally generated, but $c_4(S^*)=0$; see
\cite{ottaviani}. 
This implies that a generic section of $S^*$  does not vanishes anywhere. Choosing
such a section, we get an exact sequence
$$0\ra\cO_{\QQ_5}\ra S^*\ra Q^*\ra 0,$$
where $Q^*$ is a globally generated vector bundle and $\det Q^*=
\cO_{\QQ_5}(2)$. Moreover, there is an exact sequence 
$0\ra \we^2Q\ra \we^2S\ra Q\ra 0$.
For $i\ge 3$, we have $H^i(K_{\QQ_5}\ot Q(2))\cong H^i(K_{\QQ_5}\ot
S(2))$, which we proved to be zero. Hence 
$$H^4(\QQ_5,K_{\QQ_5}\ot \we^2S(2))\cong H^4(\QQ_5,K_{\QQ_5}\ot \we^2Q(2))
\cong H^4(\QQ_5,K_{\QQ_5}\ot Q^*).$$
And since $H^i(\QQ_5,K_{\QQ_5}\ot S^*)=H^i(\QQ_5,K_{\QQ_5}\ot S(1))=0$
for $i\ge 4$, we have
$$H^4(\QQ_5,K_{\QQ_5}\ot Q^*)\cong H^5(\QQ_5,K_{\QQ_5})=\CC,$$
as claimed. Since the argument ``pulls back'' to $Y$ by $f$, we are done.
\qed

\subsection{Quadrics of dimension six}

On $\QQ_6$ the situation is very close to that of $\QQ_5$, except that
we have now two half spinor bundles $S'$ and $S''$, both of rank four,
both having for restriction to $\QQ_5$ the spinor bundle $S$ on this 
subquadric. Since $\QQ_5$ is a hyperplane section of $\QQ_6$, the 
exact sequence $0\ra\cO_{\QQ_6}(-1)\ra\cO_{\QQ_6}\ra\cO_{\QQ_5}\ra 0$
allows to deduce almost immediately from \ref{cmq5}, a 
Castelnuovo-Mumford criterion on $\QQ_6$: a coherent sheaf 
$\cF$ on $\QQ_6$ will be generated by its global sections as soon as 
$$\begin{array}{lcll}
H^{i+j}(\QQ_6,\cF(-l))=0 & {\rm for} & (i,l)=(1,1),(2,2),(4,3),(6,4) & j\ge 0,\\
H^{i+j}(\QQ_6,\cF\ot S'(-l))=0 & {\rm for} & (i,l)=(1,0),(2,1),  & j\ge 0,
\quad {\rm and}\\
H^i(\QQ_6,\cF\ot \we^2S'(-l))=0 & {\rm for} & (i,l)=(3,1),(4,2),(5,3). & 
\end{array}$$ 

\medskip\noindent {\em Proof of Theorem C for $n=6$}.
The argument is almost exactly the same as on $\QQ_5$. Indeed, 
we have $c_4(S')=c_4(S'')=0$ \cite{ottaviani}. 
Choosing a nowhere vanishing section of 
$S'^*$, we get an exact sequence 
$$0\ra\cO_{\QQ_6}\ra S'^*\ra Q'^*\ra 0,$$
where $Q'^*$ is a globally generated vector bundle and $\det Q'^*=
\cO_{\QQ_6}(2)$. Using this sequence, we check that 
$$H^5(\QQ_6,K_{\QQ_6}\ot \we^2S'(2))=H^5(\QQ_6,K_{\QQ_6}\ot Q'^*)=
H^5(\QQ_6,K_{\QQ_6})=\CC.$$
This is enough to get the result that if $f :Y\ra\QQ_6$ is a branched covering 
with its associated vector bundle $\cE$, then $\cE(-1)$ is generated by its
global sections, 
since the other cohomology groups involved in our Castelnuovo-Mumford 
criterion on $\QQ_6$ are easily checked to vanish. \qed

\section{Global generation}

In this section we prove that for any branched covering of ordinary and
symplectic flag manifolds, the associated vector bundle $\cE$ is 
generated by its global sections. Of course, $\cE$ cannot be 
ample in general, since flag manifolds that are not grassmannians
fiber over smaller flag manifolds, from which one could pull-back 
any branched covering. Nevertheless, they should have intermediate 
positivity properties, presumably some $k$-ampleness in the sense of
Sommese. This would again imply Barth-Lefschetz type theorems, but 
unfortunately we have been unable to prove such positivity properties.

\subsection{Ordinary flag manifolds}

Let $\FF_n$ denote the variety of complete flags in
$V=\CC^n$, with universal flag
$$W_\bullet : \quad\quad 0=W_0 \subset W_1  \subset \cdots 
\subset W_{n-1} \subset V\ot \cO_{\FF_n}$$
of subbundles $W_i$ of rank $i$ of the trivial vector bundle
$V\ot \cO_{\FF_n}$. We denote by $Q_i=V/W_{n-i}$ 
the quotient bundle of rank $i$, and
$$\cO_{\FF_n}(a_1,\ldots ,a_{n-1})=\cO_{\FF_n}(\sum_{i=1}^{n-1}a_i\e_i)=
\bigotimes_{i=1}^{n-1}(\det Q_i)^{a_i}.$$
This line bundle is generated by its global sections (ample, respectively)
if and only if $a_i\ge 0$ ($a_i>0$, respectively) for all $i$.  
The canonical line bundle of $\FF_n$ is 
$K_{\FF_n}=\cO_{\FF_n}(-2,\dots,-2)$.

Let $x$ be a point of $\FF_n$, given by some complete flag 
$0=U_0 \subset U_1  \subset \cdots \subset U_{n-1} \subset V$. 
Let us choose a compatible basis $u_1,\ldots ,u_n$ of $V$, that is
such that $u_1,\ldots ,u_i$ is a basis of $U_i$ for all $i$. 
Each vector $u_{n-i}$ defines a global section of the quotient vector
bundle $Q_i$. We therefore have a global section of 
$Q_1\oplus \cdots \oplus Q_{n-1}$, vanishing precisely at $x$. 
The corresponding  Koszul complex is 
$$0\ra \we^{{n(n-1)}\over 2} (Q_1^*\oplus \cdots \oplus Q_{n-1}^*) 
\ra \cdots \ra Q_1^*\oplus \cdots \oplus Q_{n-1}^* \ra \cO_{\FF_n} 
\ra \cO_x \ra 0.$$
Hence we obtain a very simple Castelnuovo-Mumford criterion on $\FF_n$:

\begin{prop} A coherent sheaf $\cF$ on $\FF_n$ is generated by its  global
  sections as soon as 
$$H^q(\FF_n,  \we^{i_1}Q_1^*\otc \we^{i_{n-1}}Q_{n-1}^*\ot\cF)=0$$
for $q=i_1+\cdots +i_{n-1}>0$. \end{prop}

We will apply this criterion to the vector bundle associated with a
finite covering of $\FF_n$. With the generalized Le Potier vanishing theorem 
\ref{lp} in mind, it will be easy to prove the first half of Theorem D:

\begin{prop} Let $f: Y\ra \FF_n$ be a finite surjective morphism 
where $Y$ is a nonsingular connected complex projective variety.
Then the associated vector bundle $\cE$ is generated by its global
sections. \end{prop}

\noindent {\em Proof.} We need to prove that 
$$H^q(\FF_n, K_{\FF_n} \ot
\we^{i_1} Q_1^* \otc \we^{i_{n-1}} Q_{n-1}^* (2,\ldots ,2)\ot\cE)=0$$
for $q=i_1+\cdots +i_{n-1}>0$.
Using Serre duality, the definition of $\cE$, and Leray spectral sequence, 
this vanishing is 
seen to be equivalent to the equality between the cohomology group
$$H^q(\FF_n, K_{\FF_n} \ot
\we^{i_1} Q_1^* \otc \we^{i_{n-1}} Q_{n-1}^* (2,\ldots ,2)),$$
and the corresponding cohomology group on $Y$, 
$$H^q(Y, K_Y \ot
f^*(\we^{i_1} Q_1^* \otc \we^{i_{n-1}} Q_{n-1}^* (2,\ldots ,2)))$$
for $q=i_1+\cdots +i_{n-1}>0$.
We shall prove that both groups are equal to zero. We begin with 
the first one. 

If $i_1 >0$, then since $Q_1$ is of rank $1$ we can suppose that
$i_1=1$, and we apply the generalized Le Potier vanishing 
theorem \ref{lp}. 
If $i_1=0$, let $k$ be the smallest integer such that $i_k >0$. 
We proceed by induction on $k$. We have a short exact sequence
$$0\ra \cO_{\FF_n}(\e_k-\e_{k-1}) \ra Q_k \ra Q_{k-1}\ra 0,$$
from which we get an exact sequence for wedge powers:
$$0\ra \we^{i_k} Q_{k-1}^* \ra \we^{i_k} Q_k^* \ra 
\we^{i_k -1} Q_{k-1}^*(\e_{k-1}-\e_k )\ra 0.$$
Let us tensor this exact sequence  by 
 $K_{\FF_n}\ot \we^{i_{k+1}} Q_1^* \otc
\we^{i_{n-1}} Q_{n-1}^*(2,\dots ,2)$. The $q$-th cohomology group of 
the first vector bundle of the exact sequence we get, then vanishes 
by induction hypothesis.
Also, for the last vector bundle we obtain the cohomology group 
$$H^q(\FF_n, K_{\FF_n}\ot\we^{i_{k-1}} Q_{k-1}^* \ot \we^{i_{k+1}} Q_{k+1}^*
\cdots \ot \we^{i_{n-1}} Q_{n-1}^* (2,\dots,3,1,2,\dots, ,2)),$$
where $3$ and $1$ occupy the $(k-1)$-th and $k$-th position, respectively.
But this is zero by the generalized Le Potier vanishing theorem. 
The $q$-th cohomology group of the middle vector bundle therefore 
vanishes, which is what we wanted to prove. 

Since the very same argument applies to the pulled-back vector bundles
on  $Y$, the proposition is proved. \qed

\begin{coro}\label{flagA} Let $f: Y\ra \FF$ be a finite covering of any variety 
of incomplete flags $\FF$, where $Y$ is a nonsingular connected complex
 projective variety.
Then the associated vector bundle $\cE$ is generated by its global
sections. \end{coro}

\noindent {\em Proof}. Since a variety of complete flags fibers on 
any variety of incomplete flags of the same vector space, we have a 
smooth fibration $\pi :\FF_n\ra\FF$. Pulling back by $f$, we get 
 a finite covering $f': Y'\ra\FF_n$ with its  associated 
vector bundle  $\cE'=\pi^*\cE$. Moreover, since $\pi$ has
connected fibers, $\pi_*\cO_{\FF_n}=\cO_{\FF}$, hence 
$$H^0(Y',\cE')=H^0(Y,\pi_*\cE')=H^0(Y,\cE).$$
Since, by the proposition above, $\cE'$ is generated by its global
sections, all of which are pull-backs of sections of $\cE$, 
 then $\cE$ itself must be globally generated. \qed

\subsection{Symplectic flag manifolds}

We denote by ${\LF}_n$ the variety of complete 
lagrangian flags in $V=\CC^{2n}$, with universal flag of subbundles 
$$W_\bullet : \quad 0=W_0 \subset W_1  \subset \cdots \subset 
W_n = W_n^\perp\subset W_{n-1}^\perp\subset \cdots \subset 
W_1^\perp\subset \cO_{\LF_n}\ot V,$$
where $\dim W_i = i$ and $\dim W_i^{\perp} = 2n-i$. 
This variety has dimension $n^2$,
and is homogeneous under the natural action of the symplectic group
$Sp(2n,\CC)$. We will consider it as a subvariety of $\FF_{n,2n}$, the
variety of incomplete flags in $\CC^{2n}$ consisting of subspaces of
dimensions $1$ to $n$. 

For $1\le i\le n$, we denote by $Q_i=V/W_i$
the quotient bundle of rank $2n-i$, and
$$\cO_{\LF_n}(a_1,\ldots ,a_n)=\cO_{\FF_n}(\sum_ia_i\e_i)=
\bigotimes_{i=1}^n(\det Q_i)^{a_i}.$$
In particular, the canonical line bundle of $\LF_n$ is 
$K_{\LF_n}=\cO_{\LF_n}(-2,\dots,-2)$.

On $\FF_{n,2n}$, a generic section of $Q_1\oplus\cdots\oplus Q_n$
vanishes at a simple point, and the associated Koszul complex provides
us with a Castelnuovo-Mumford criterion, which  can be restricted to 
$\LF_n$. Since  we have $Q_i^*\simeq W_{n-i+1}^{\perp}$ on $\LF_n$, we
get:

\begin{prop} A coherent sheaf $\cF$ on $\LF_n$ is generated by its global
  sections as soon as 
$$H^q(\LF_n,  \we^{i_1}W_1^{\perp}\otc \we^{i_n}W_n^{\perp}\ot\cF)=0$$
for $q=i_1+\cdots +i_n>0$. \end{prop}

In the sequel, we will need to understand the cohomology of certain
homogeneous bundles, which will be tensor products of wedge powers 
of the tautological vector bundles. We  call such a bundle 
a {\em wedge bundle}, and  say that it is of 
{\em type} $t=(t_1,\ldots ,t_n)$ if it has $t_i$ factors that are 
wedge powers of  $W_i$ or $W_i^{\perp}$, and {\em weight} $w=(w_1,\ldots,w_n)$
if the sum of the exponents of these wedge powers of $W_i$ and $W_i^{\perp}$
is equal to $w_i$. The type can of course be decomposed into subtypes
$t'$ and $t''$,  given by the number of factors of involving
$W_i$, and $W_i^{\perp}$, respectively. 
We  first need the following refined statement of the generalized Le Potier
vanishing theorems for wedge bundles on $\LF_n$.

\begin{prop}\label{argh} Let $\cW$ be a wedge bundle on $\LF_n$, of type 
$t$ and weight $w$. Then 
$$H^q(\LF_n,K_{\LF_n}\ot \cW(l_1,\ldots ,l_n))=0$$
for $q>w_1+\cdots +w_n$, if either :
\begin{itemize}
\item $l_k>t_k$ for $1\leq k\leq n$,
\item $l_k>t_k$ for $k\ne p,q$, where $p<q$ are such that 
$t''_q>0$, $l_q\ge t_q$ and $l_p>t_p+1,$ or
\item $l_k>t_k$ for $k\ne p,q$, where $p>q$ are such that 
$t'_q>0$, $l_q\ge t_q$ and $l_p>t_p+1$.
\end{itemize}\end{prop}

\noindent {\em Proof}. The first assertion is an immediate consequence
of the generalized Le Potier vanishing theorem. Indeed, we can rewrite
each wedge power $\we^mW_i$ as $\we^{n-i+1-m}W_i^*(-1)$, and similarly
each $\we^mW_i^{\perp}$ as $\we^{n-i+1-m}(W_i^{\perp})^*(-1)$. If
$l_k>t_k$ for $1\leq k\leq n$, we get a tensor product of wedge 
powers of globally generated vector bundles, by the ample line bundle
$\cO_{\LF_n}(l_1-t_1,\ldots ,l_n-t_n)$. 
Hence we can apply the generalized Le Potier vanishing theorem
\ref{lp}, and we obtain the desired statement. 

For the second assertion, which is not a direct consequence of Le Potier's
theorem, we will proceed by induction on $q$. Since $t''_q$ is non zero, 
then we
have at least a factor $\we^mW^{\perp}_q$ in $\cW$. Therefore we can use 
the exact sequence 
$$0\ra\we^m W_q^{\perp}\ra\we^m W_{q-1}^{\perp}\ra
\we^{m-1} W_q^{\perp}(\e_q-\e_{q-1})\ra 0.$$
Let us tensor this exact sequence by the other factors of $\cW$. 
We get a short exact sequence involving wedge bundles,
$$0\ra\cW\ra\cW'\ra\cW''(\e_q-\e_{q-1})\ra 0.$$
Here $\cW'$ has weight $(w_1,\ldots ,w_{q-1}+m,w_q-m,\ldots ,w_n)$ 
and type $(t_1,\ldots ,t_{q-1}+1,t_q-1,\ldots ,w_n)$. Hence 
$$H^q(\LF_n,K_{\LF_n}\ot \cW'(l_1,\ldots ,l_n))=0$$
for $q>w_1+\cdots +w_n$. Indeed, this is the case by induction
hypothesis if $p<q-1$, while if $p=q-1$ we are in the situation of the
first part of the proposition. On the other hand, 
$\cW''$ has weight $(w_1,\ldots ,w_{q-1},w_q-1,\ldots ,w_n)$ 
and type $(t_1,\ldots ,t_{q-1},t_q,\ldots ,w_n)$ (except if $m=1$, in
which case $t_q$ is replaced by $t_q-1$, which is even better). This 
again allows us to use our induction hypothesis, which gives 
$$H^{q-1}(\LF_n,K_{\LF_n}\ot \cW''(l_1,\ldots ,l_{q-1}-1,
l_q+1,\ldots ,l_n))=0$$
for $q>w_1+\cdots +w_n$. And these vanishings for $\cW'$ and $\cW''$
imply our claim for $\cW$. 

The third assertion is proved in the very same way using the other
series of exact sequences,
$$0\ra\we^m W_q\ra\we^m W_{q+1}\ra
\we^{m-1} W_q(\e_q-\e_{q+1})\ra 0,$$
and a similar, but now  descending, induction on $q$.\qed

\medskip With the help of this vanishing theorem, we can now extend
\ref{flagA} to the case of symplectic flag manifolds and prove the
second half of Theorem D.
As for ordinary flags, it will suffice to consider complete symplectic
flags, that is, coverings of $\LF_n$. 

\begin{prop} Let $f: Y\ra\LF_n$ be a finite covering, where $Y$ is a
 nonsingular connected complex projective variety. Then the associated vector bundle $\cE$ is 
  generated by its global sections. \end{prop}

\noindent {\em Proof}. We apply our Castelnuovo-Mumford criterion on
  $\LF_n$ to $\cE$. Because of Serre duality and the definition of $\cE$, 
 we need to prove  that the cohomology group 
$$H^q(\LF_n, K_{\LF_n} \ot
\we^{i_1} W_1^{\perp} \otc \we^{i_n} W_n^{\perp} (2,\ldots ,2)),$$
for $q=i_1+\cdots +i_n>0$,
is equal to the corresponding group obtained  by pull-back to $Y$. 
We shall prove that both groups are equal to zero. Note that on $\LF_n$,
this would follow from the first part of the proceeding proposition if
we had $W_1$ instead of $W_1^{\perp}$, at least for $i_1>0$. Indeed,
since $W_1$ has rank one, this would force $i_1=1$ and we could
replace $W_1$ by $\cO_{\LF_n}(-\e_1)$. Our strategy will therefore be to use
the tautological exact sequences on $\LF_n$, first to replace
$W_1^{\perp}$ by $W_n^{\perp}=W_n$, then $W_n$ by $W_1$. 

Let $k$ be the smallest integer such that 
$i_k>0$. For notational simplicity, let us suppose that $k=1$: 
the argument would be the same for $k>1$. We have an exact sequence
$$0\ra\we^{i_1} W_2^{\perp}\ra\we^{i_1} W_1^{\perp}\ra
\we^{i_1-1} W_2^{\perp}(\e_2-\e_1)\ra 0.$$
Let us tensor it by $\cW=\we^{i_2} W_2^{\perp} \otc \we^{i_n} 
W_n^{\perp}$. The wedge bundle $\cW\ot\we^{i_1-1}
W_2^{\perp}(\e_2-\e_1)$ has type $(0,2,1,...,1)$; so that 
by the second part of proposition \ref{argh},
$$H^q(\LF_n,K_{\LF_n}\ot \cW\ot\we^{i_1-1}
W_2^{\perp}(1,3,2,\ldots ,2)=0.$$
Our vanishing will therefore follow from that of
$H^q(\LF_n, K_{\LF_n} \ot \we^{i_1} W_2^{\perp}\ot\cW (2,\ldots ,2))$.
Using the same argument $n-1$ times, we are then reduced to prove that
$$H^q(\LF_n, K_{\LF_n} \ot \we^{i_1} W_n^{\perp}\ot\cW(2,\ldots ,2))=0.$$
But since $W_n^{\perp}=W_n$, we can now use the exact sequence
$$0\ra\we^{i_1} W_{n-1}\ra\we^{i_1} W_n\ra
\we^{i_1-1} W_{n-1}(\e_{n-1}-\e_n)\ra 0.$$
Using the third part of proposition \ref{argh}, we see that 
$$H^q(\LF_n, K_{\LF_n} \ot \we^{i_1} W_{n-1}\ot\cW(2,\ldots ,2,3,1))=0.$$
We are thus reduced to prove that
$H^q(\LF_n, K_{\LF_n} \ot \we^{i_1} W_{n-1}\ot\cW(2,\ldots ,2))=0$
and, using repeatedly the same argument, what we finally have
to show is the vanishing of
$$H^q(\LF_n, K_{\LF_n} \ot \we^{i_1} W_1\ot\cW(2,\ldots ,2)).$$
This is clear if $i_1>1$. If $i_1=1$, we can rewrite this cohomology
group as
$$H^q(\LF_n, K_{\LF_n} \ot\cW(1,2,\ldots ,2)),$$
which is zero by the first part of \ref{argh}. \qed

\subsection{Questions and conjectures}

Let $X=G/P$ be a complex projective homogeneous space, where $G$ is a complex
semisimple Lie group and $P$ is a parabolic subgroup. If $P$ is not
maximal and $Q$ is a parabolic subgroup of $G$ containing $P$, we
have a smooth fibration $G/P\ra G/Q$. Then we obtain a branched covering  $f$
of $X$  by pulling back the one of $G/Q$. Note that its associated 
bundle, being a pull-back by $f$, cannot be ample: it is at best 
$(\dim Q/P)$-ample. Let us
define $$k_X=\max_{Q\supset P}\dim Q/P.$$
If $P$ is defined (up to conjugacy) by a set $I$ of simple roots, 
one can compute $k_X$ as the maximum, as $\a$ describes $I$, of the 
number of positive roots having positive coefficient on the simple 
root $\a$, but zero coefficient on the other simple roots in $I$. 

We can then extend the conjecture of Debarre stated in the introduction 
in the following naive way:

\medskip\noindent {\bf Conjecture}. {\em Let $f :Y\ra X=G/P$ be a
  branched covering, with $Y$ smooth and connected. Then the
  associated vector bundle $\cE$ is $k_X$-ample.}

\medskip\noindent Once again, this would imply a Barth-Lefschetz
type theorem: the natural map 
$$f^* : H^i(X,\CC)\lra H^i(Y,\CC)$$
would then be an isomorphism for $i\leq \dim X -k_X-d+1$, 
where $d$ is the degree of $f$. 
Moreover, one could expect this $k_X$-ampleness to be optimal only 
when $f$ is obtained by pull-back:

\medskip\noindent {\bf Question}. {\em Let $f :Y\ra X=G/P$ be a
  branched covering, with a nonsingular connected complex projective variety
$Y$. Suppose that the
  associated vector bundle $\cE$ is not $(k_X-1)$-ample. Then, is 
$f$ necessarily  a pull-back of a covering of some smaller homogeneous
  space $G/Q$ ?}

\medskip To attack the previous conjecture by the methods of this
paper, we would need {\em efficient} Castelnuovo-Mumford criteria. Even 
when $P$ is maximal, we could not always find good enough such criteria,
for example on quadrics, or on spinor varieties, which being hermitian
symmetric should be easier to deal with than orthogonal grassmannians. 

\medskip\noindent {\bf Question}. {\em How to find ``good'' 
Castelnuovo-Mumford criteria on homogeneous spaces ?}

\medskip\noindent
When the parabolic group $P$ is not maximal, such criteria should 
certainly involve several line bundles. A special case would be the 
following:

\medskip\noindent {\bf Question}. {\em Let $L,M$ be globally generated 
line bundles on a projective variety $X$, with $L\ot M$ very ample. 
If $\cF$ is a coherent sheaf on $X$, can one give a cohomological 
criterion for the surjectivity of the mixed evaluation morphism
$$H^0(X,\cF\ot L^{-1})\ot L\;\op \; H^0(X,\cF\ot M^{-1})\ot M\lra \cF,$$
that would not necessarily imply that $\cF\ot L^{-1}$ or 
$\cF\ot M^{-1}$ is globally generated ?}

\noindent\begin{tabular}{ll}
 Meeyoung Kim&Laurent Manivel\\
Max-Planck-Institut f\"ur Mathematik\quad &Institut Fourier, UMR 5582 UJF/CNRS\\
Gottfried-Claren-Str. 26&Universit\'e de Grenoble I, BP 74 \\
D-53225 Bonn,
Germany\ \ \ \ \ \ \ \ \ &F-38402 Saint Martin d'H\`eres, France\\
kim@mpim-bonn.mpg.de&laurent.manivel@ujf-grenoble.fr\\
\end{tabular}

\end{document}